\documentclass[12pt]{article}
\usepackage{amssymb}
\usepackage{amsopn}
\usepackage{latexsym}
\usepackage{esint}
\usepackage{amsfonts}
\title{\textbf{A first regularity result for the Armstrong-Frederick cyclic hardening plasticity model with Cosserat effects}}
\author{\textbf{Krzysztof Che{\l}mi\'{n}ski}\\[0.5ex]
\textbf{\footnotesize{Faculty of Mathematics and
 Information Science, Warsaw University of Technology,}}\\[-1ex]
\textbf{\footnotesize{ul. Koszykowa 75, 00-662 Warsaw, Poland}}\\[-1ex]
\textbf{\footnotesize{E-Mail: kchelmin@mini.pw.edu.pl}}\\[1ex]
\textbf{Patrizio Neff}\\[0.5ex]
\textbf{\footnotesize{Head of Lehrstuhl f\"ur Nichtlineare Analysis und Modellierung,}}\\[-1ex] 
\textbf{\footnotesize{Fakult\"at f\"ur Mathematik, Universit\"at
    Duisburg-Essen}}\\[-1ex]
\textbf{\footnotesize{Thea-Leymann Strasse 9,  45141 Essen, Germany}}\\[-1ex]
\textbf{\footnotesize{E-Mail: patrizio.neff@uni-due.de}}\\[1ex]
\textbf{Sebastian Owczarek}\\[0.5ex] 
\textbf{\footnotesize{Faculty of Mathematics and
 Information Science, Warsaw University of Technology,}}\\[-1ex]
\textbf{\footnotesize{ul. Koszykowa 75, 00-662 Warsaw, Poland}}\\[-1ex]
\textbf{\footnotesize{E-Mail: s.owczarek@mini.pw.edu.pl}}}
\date{}
\newtheorem{tw}{Theorem}[section]
\newtheorem{lem}[tw]{Lemma}
\newtheorem{de}[tw]{Definition}
\DeclareMathOperator{\dev}{dev}
\begin{document}
\maketitle
\begin{abstract}
\noindent
The purpose of this article is to prove the H\"older continuity up to the boundary of the displacement vector and the microrotation matrix for the quasistatic, rate-independent Armstrong-Frederick cyclic hardening plasticity model with Cosserat effects. This model is of non-monotone and non-associated type. In the case of two space dimensions we use the hole-filling technique of Widman and the Morrey's Dirichlet growth theorem. 
\end{abstract}
\newcommand{\bl}{\backslash}
\newcommand{\nn}{\nonumber}
\newcommand{\KK}{\sigma_{\rm y}}
\newcommand{\ve}{\varepsilon}
\newcommand{\R}{{\mathbb R}}
\newcommand{\D}{{\mathbb C}}
\newcommand{\E}{{\cal E}}
\newcommand{\K}{{\cal K}}
\newcommand{\di}{{\mathrm d}}
\newcommand{\so}{\mathfrak{so}}
\renewcommand{\S}{{\cal S}^3}
\renewcommand{\SS}{{\cal S}^3_{\mathrm{dev}}}
\newcommand{\id}{ {1\!\!\!\:1 } }

\section{Introduction}
In the paper we investigate the important regularity question for models elasto-plasticity. Various model systems in the finite strain and small strain case have been proposed in the articles \cite{28}-\cite{31}, \cite{33}. In this contribution we focus on the infinitesimal plasticity models in the framework proposed by D. H. Alber and his group (see \cite{2}, \cite{34}, \cite{35}, \cite{36}, \cite{8}). This framework is perfectly adapted for inelastic deformation processes of metals that are characterized by monotone flow rule (associated plasticity). In that case the finite difference method was very useful to prove regularity of stresses in the Prandtl-Reuss and Norton-Hoff models \cite{4,5,13,16,26}, because this method allows to cancel the monotone nonlinearities. Using this method D. H. Alber and S. Nesenenko \cite{3} have shown $L^{\infty}(H^{1/3-\delta})$- regularity for stresses and plastic strain for coercive models of viscoplasticity with variable coefficients. Next, D. Knees in \cite{22} obtains the stresses in the space $L^{\infty}(H^{1/2-\delta})$. A similar result was proved by P. Kami\'nski in \cite{20} and \cite{21} for coercive and self-controlling (non-coercive) viscoplastic models. Moreover, in \cite{10} a $H^1_{loc}$- regularity result for the stresses and strains in Cosserat elasto-plasticity was proved, cf. \cite{25}. See also \cite {37} and \cite {38} for the local regularity in the Hencky model. 

The Armstrong-Frederick (AF) cyclic hardening plasticity model with {C}osserat effects was formulated first in the article \cite{12} and it is of non-monotone and not of gradient type (non-associated flow rule - see \cite{7} and \cite{12} for more details). Hence the finite difference method is useless to study regularity of solutions for this model. The idea of the present paper is to show the H\"older continuity up to the boundary of the displacement and the microrotation matrix in AF-model with Cosserat effects in the case of two space dimensions (the existence of the energy solutions for this model was also proved in \cite{12}).  We will use a very old method, which was exposed in \cite{23} and \cite{24}. Those works presented first fundamental theorems about existence and regularity of solutions of two - dimensional elliptic systems.

We derive the Morrey's condition up to the boundary of the basic domain $\Omega$. In order to got it we will use the Widman's hole filling trick proposed by K. O. Widman in the paper \cite{27}.

Morrey's methods were used in \cite{15} where the author shown the existence of a H\"older continuous solution for a class of two-dimensional non-linear elliptic systems. Moreover this methods were also used to prove H\"older continuity for the displacements in isotropic and kinematic hardening with von Mises yield criterion in \cite{17}.

To our knowledge this article presents the first regularity result for models from the theory of inelastic deformations of metals, which are non-monotone, non-associated but coercive - for the definitions we refer to \cite{2}.  

%%%%%%%%%%%%%%%%%%%%%%%%%%%%%%%%%%%%%%%%%%%%%%%%%%%%% Section 2 %%%%%%%%%%%%%%%%%%%%%%%%%%%%%%%%%%%%%%%%%%%%%%%%%%%%%%%%%%%%%%%%%%%%%%%%

\section{The Armstrong-Frederick model with Cosserat effect}
Here we formulate the Armstrong-Frederick model with Cosserat effects in the case of two space dimensions.

Let us assume that $\Omega\subset\R^2$ is a bounded domain with Lipschitz boundary $\partial\Omega$. The structure of the model equations is the same as introduced in \cite{12}. For the mechanical results for Cosserat plasticity we refer to \cite{8}, \cite{9} - see also \cite{11}, where the non-monoton model of poroplastisity with Cosserat effects was introduced. Hence we deal with the following system of equations  
\renewcommand{\theequation}{\thesection.\arabic{equation}}
\setcounter{equation}{0}%
\begin{eqnarray}
\label{eq:2.1}
\mathrm{div}_x\, T&=&-f\,,\nn\\
T&=&2\mu(\ve(u)-\ve^p)+2\mu_c(\mathrm{skew}\,(\nabla_x u)-A)+\lambda\mathrm{tr}\,(\ve(u)-\ve^p)\id\,,\nn\\[1ex]
-l_c\,\Delta_x\,\mathrm{axl}\,(A)&=&\mu_c\,\mathrm{axl}\,(\mathrm{skew}\,(\nabla_x u)-A)\,,\nn\\[1ex]
\ve^{p}_{t}&\in&\partial I_{K(b)}\Big (T_{E}\Big)\,,\\
T_E&=&2\mu(\ve(u)-\ve^p)+\lambda\mathrm{tr}\,(\ve(u)-\ve^p)\id,\nn\\
b_t&=&c\,\ve^p_t-d\,|\ve_t^p|b\,,\nn
\end{eqnarray}
where the unknowns are: the displacement vector field $u:\Omega\times[0,T]\rightarrow \R^3$, the microrotation matrix
$A:\Omega\times[0,T]\rightarrow\mathfrak{so}(3)$ ($\mathfrak{so}(3)$ is the set of skew-symmetric $3\times 3$ matrices) and the vector of internal variables $z=(\ve^p,b):\Omega\times [0,T]\rightarrow\SS\times\SS$ ($\ve^p$ is the classical infinitesimal symmetric plastic strain tensor, $b$ is the symmetric backstress tensor and the space $\SS$ denotes the set of symmetric $3\times 3$ matrices with vanishing trace). $\ve(u)=\mathrm{sym}\,(\nabla_x u)$ denotes the symmetric part of the gradient of the displacement.

The equations (\ref{eq:2.1}) are studied for $x\in\Omega\subset \R^2$ and $t\in [0,T]$, where $t$ denotes the time.

The set of admissible elastic stresses $K(b(x,t))$ is defined in the form\\
$K(b)=\{T_E\in\S\,:\,|\dev\,(T_E)-b|\leq \KK\}$, where $\dev\,(T_E)=T_E-\frac{1}{3}\,\mathrm{tr}\,(T_E)\cdot\id$, $\KK$ is a material parameter (the yield limit) and $\id$ denotes the identity $3\times 3$ matrix. The function $I_{K(b)}$ is the indicator function of the admissible set $K(b)$ and $\partial I_{K(b)}$ is the subgradient of the convex, proper, lower semicontinous function $I_{K(b)}$.

The function $f:\Omega\times[0,T]\rightarrow \R^3$ describes the density of the applied body forces, the parameters $\mu$, $\lambda$ are positive Lam\'e constants, $\mu_c>0$ is the Cosserat couple modulus and $l_c>0$ is a material parameter with dimensions $[m^2]$, describing a length scale of the model due to the Cosserat effects.  $c,d>0$ are material constants.

The operator $\mathrm{skew}\,(T)=\frac{1}{2}(T-T^T)$ denotes the skew-symmetric part of a $3\times 3$ tensor. The operator $\mathrm{axl}:\mathfrak{so}(3)\rightarrow \R^3$ establishes the identification of a skew-symmetric matrix with a vector in $\R^3$. This means that if we take $A\in\mathfrak{so}(3)$, which is in the form $A=((0,\alpha,\beta),(-\alpha,0,\gamma),(-\beta,-\gamma,0))$, then $\mathrm{axl}(A)=(\alpha,\beta,\gamma)$.\\[1ex]
Notice that the system (\ref{eq:2.1}) is a modification of the Melan-Prager model, which is well known in the literature and it can also be seen as an approximation of the Prandtl-Reuss model. The expression $|\ve^p_t|b$ is a perturbation of the Melan-Prager model - if $d=0$ then we obtain the classical Melan-Prager linear kinematic hardening model - details can be found in \cite{7} and \cite{12}.\\[1ex]
The system (\ref{eq:2.1}) is considered with the Dirichlet boundary condition for the displacement:
\begin{eqnarray}
\label{eq:2.2}
u(x,t)=g_D(x,t)\qquad \textrm{ for}\quad x\in\partial\Omega\quad\textrm{and}\quad t\geq 0
\end{eqnarray}
and with the Dirichlet boundary condition for the microrotation:
\begin{eqnarray}
\label{eq:2.3}
A(x,t)=A_D(x,t)\qquad \textrm{ for}\quad x\in \partial \Omega \quad\textrm{and}\quad  t\geq 0.
\end{eqnarray}
Finally, we consider the system (\ref{eq:2.1}) with the following initial conditions
\begin{eqnarray}
\label{eq:2.4}
\ve^p(x,0)=\ve^{p,0}(x),\qquad b(x,0)=b^0(x).
\end{eqnarray}
The free energy function associated with the system (\ref{eq:2.1}) is given by the formula
\begin{eqnarray}
\label{eq:2.5}
\rho\,\psi(\ve,\ve^p,A,b)&=&\mu\,\|\ve(u)-\ve^p\|^2+
\mu_c\,\|\mathrm{skew}\,(\nabla_x u)-A\|^2\nn\\[1ex]
&+&\frac{\lambda}{2}\,\Big(\mathrm{tr}\,(\ve(u)-\ve^p)\Big)^2+2\,l_c\,\|\nabla_x\,\mathrm{axl}\,(A)\|^2+\frac{1}{2c}\,\|b\|^2,
\end{eqnarray}
where $\rho$ is the mass density which we assume to be constant in time and space. The total energy is of the form:
\begin{eqnarray*}
\E(\ve,\ve^p,A,b)(t)=\int\nolimits_{\Omega}\rho\psi(\ve(x,t),\ve^p(x,t),A(x,t),b(x,t))\,\di x.
\end{eqnarray*}
The Section $2$ of the article \cite{7} shows that the inelastic constitutive equation occurring in (\ref{eq:2.1}) is of pre-monotone type (for the definition see \cite{2}). We also know that the AF-model with micropolar effects is of non-monotone type and not of gradient type (non-associated flow rule).

%%%%%%%%%%%%%%%%%%%%%%%%%%%%%%%%%%%%%%%%%%%%%%%%%%%%% Section 3 %%%%%%%%%%%%%%%%%%%%%%%%%%%%%%%%%%%%%%%%%%%%%%%%%%%%%%%%%%%%%%%%%%%%%%%%

\section{Existence theory and main result}
The only one existence result for the AF-models with Coserrat effects in the case of three space dimensions was obtained in the article \cite{12}. It was shown that the limit in the Yosida approximation process satisfies the energy inequality for special test functions. Using the same techniques as in \cite{12} we can obtain the same existence theorem in the two dimensional cases. The goal of this section is to formulate the existence theorem for the system (\ref{eq:2.1}) and the main result of this article.\\
\renewcommand{\theequation}{\thesection.\arabic{equation}}
\setcounter{equation}{0}%
Let us  assume that for all $T>0$ the given data $f$, $g_D$ and $A_D$ have the regularity
\begin{equation}
\label{eq:3.1}
f\in H^{1}(0,T;L^2(\Omega;\R^3)), \quad g_D\in H^{1}(0,T;H^{\frac{1}{2}}(\partial\Omega;\R^3))\,,
\end{equation}
\begin{equation}
\label{eq:3.2}
A_D\in H^{1}(0,T;H^{\frac{3}{2}}(\partial\Omega;\so(3))).
\end{equation}
Additionally let us suppose that the initial data $(\ve^{p,0},b^0)\in L^2(\Omega;\SS)\times L^2(\Omega;\SS)$\\ satisfy
\begin{equation}
\label{eq:3.3}
|b^{0}(x)|\leq\frac{c}{d}\quad \mathrm{and}\quad |\dev\,(T^0_E(x))-b^0(x)|\leq \KK \quad \mathrm{for\; almost\; all\;}x\in\Omega,
\end{equation}
where the initial stress $T^0_E=2\mu\,(\ve(u(0))-\ve^{p,0})+\lambda\,\mathrm{tr}\,(\ve(u(0))-\ve^{p,0})\,\id\in L^2(\Omega;\S)$ is the unique solution of the following linear problem
\begin{eqnarray}
\label{eq:3.4}
\mathrm{div}_x\, T^0(x)&=&-f(x,0),\nn\\[1ex]
-l_c\,\Delta_x\,\mathrm{axl}\,(A(x,0))&=&\mu_c\,\mathrm{axl}\,(\mathrm{skew}\,(\nabla_x u(x,0))-A(x,0)),\nn\\[1ex]
u(x,0)_{|_{\partial\Omega}}=g_D(x,0) && A(x,0)_{|_{\partial\Omega}}=A_D(x,0),
\end{eqnarray}
with
\begin{eqnarray*}
T^0(x)=2\mu(\ve(u(x,0))-\ve^{p,0}(x))&+&2\mu_c(\mathrm{skew}\,(\nabla_x u(x,0))-A(x,0))\\[1ex]
&+&\lambda\mathrm{tr}\,(\ve(u(x,0))-\ve^{p,0}(x))\,\id.
\end{eqnarray*}
Let us consider the convex set (which will be used as set of test functions further on)
$$\K^{\ast}=\{(\dev\,(T_E),-\frac{1}{c}b)\in\SS\times\SS\,:\, |\dev\,(T_E)-b|+\frac{d}{2c}\,|b|^2\leq \KK\},$$
where the constant $\KK$ is the same as in the yield condition. The theory of elasticity implies that there exists a positive definite operator $\D^{-1}:\S\rightarrow\S$ such that $\D^{-1}T_{E,t}=\ve_t-\ve^p_t$. Now we recall from \cite{12} a notion of the definition of the energy solution for the system (\ref{eq:2.1}) (for a motivation  we refer to \cite{12}).
\begin{de}
\label{de:3.1}
$\mathrm{(solution\; concept-energy\; inequality)}$\\
Fix $T>0$. Suppose that the given data satisfy (\ref{eq:3.1}) - (\ref{eq:3.4}). We say that a vector $(u,T,A,\ve^p,b)\in L^{\infty}(0,T;H^{1}(\Omega;\R^3)\times L^2(\Omega;\S)\times H^2(\Omega;\mathfrak{so}(3))\times (L^{2}(\Omega;\SS))^2)$ solves the problem (\ref{eq:2.1})-(\ref{eq:2.4}) if
$$(u_t,T_t,A_t,\ve^p_t,b_t)\in L^{2}(0,T;H^{1}(\Omega;\R^3)\times L^2(\Omega;\R^9)\times H^2(\Omega;\mathfrak{so}(3))\times (L^{2}(\Omega;\SS))^2),$$
the equations $(\ref{eq:2.1})_1$ and $(\ref{eq:2.1})_3$ are satisfied pointwise almost everywhere on $\Omega\times (0,T)$ and for all test functions $(\hat{T}_E,\hat{b})\in L^{2}(0,T;L^{2}(\Omega;\S)\times L^{2}(\Omega;\SS))$ such that 
$$ (\dev\,(\hat{T}_E),\hat{b})\in \K^{\ast},\qquad\qquad \mathrm{div}\; \hat{T}_E\in L^{2}(0,T;L^2(\Omega,\R^3)),$$
the inequality
\begin{eqnarray}
&&\frac{1}{2}\int\nolimits_{\Omega}\D^{-1}\,T_{E}(x,t)\,T_E(x,t)\,\di x+ \mu_c\int\nolimits_{\Omega}|\mathrm{skew}\,(\nabla_x u(x,t))-A(x,t)|^2\di x\nn\\[1ex]
&+&2\,l_c\int\nolimits_{\Omega}|\nabla\,\mathrm{axl}\,(A(x,t))|^2\,\di x+ \frac{1}{2c}\int\nolimits_{\Omega}|b(x,t)|^2\,\di x
\leq\frac{1}{2}\int\nolimits_{\Omega}\D^{-1}\,T^0_{E}\,(x)T^0_E(x)\,\di x\nn\\[1ex]
&+& \mu_c\int\nolimits_{\Omega}|\mathrm{skew}\,(\nabla_x u(x,0))-A(x,0)|^2\,\di x+ \frac{1}{2c}\int\nolimits_{\Omega}|b(x,0)|^2\,\di x\nn\\[1ex]
&+&2\,l_c\int\nolimits_{\Omega}|\nabla\,\mathrm{axl}\,(A(x,0))|^2\,\di x
+\int\nolimits_0^t\int\nolimits_{\Omega}u_t(x,\tau)\,f(x,\tau)\,\di x\,\di\tau\nn\\[1ex]
&+&\int\nolimits_0^t\int\nolimits_{\Omega}u_t(x,\tau)\,\mathrm{div}\,\hat{T}_E(x,\tau)\,\di x\,\di\tau
+\int\nolimits_0^t\int\nolimits_{\partial\Omega}g_{D,t}(x,\tau)\,(T(x,\tau)-\hat{T}_E(x,\tau))\cdot n(x)\,\di S\,\di\tau\nn\\[1ex]
&+&\int\nolimits_0^t\int\nolimits_{\Omega}\D^{-1}\,T_{E,t}(x,\tau)\,\hat{T}_E(x,\tau)\,\di x\,\di\tau
+\frac{1}{c}\int\nolimits_0^t\int\nolimits_{\Omega}b_t(x,\tau)\,\hat{b}(x,\tau)\,\di x\,\di\tau\nn\\[1ex]
&+&4\,l_c\int\nolimits_0^t\int\nolimits_{\partial\Omega}\nabla\mathrm{axl}\,(A(x,\tau))\cdot n \;\mathrm{axl}\,(A_{D,t}(x,\tau))\,\di S\,\di\tau\nn
\end{eqnarray}
is satisfied for all $t\in(0,T)$, where $T^0_E\in L^2(\Omega;\S)$ and $(u(0),A(0))\in H^{1}(\Omega;\R^3)\times H^2(\Omega;\mathfrak{so}(3))$ are unique solution of the problem (\ref{eq:3.4}).
\end{de}
\begin{tw}
$\mathrm{(Existence\; result)}$\\
\label{tw:3.2}
Let us assume that the given data and initial data satisfy the properties, which are specified in (\ref{eq:3.1}) - (\ref{eq:3.4}). Then there exists a global in time solution (in the sense of Definition \ref{de:3.1}) of the system (\ref{eq:2.1}) with boundary conditions (\ref{eq:2.2}), (\ref{eq:2.3}) and initial condition (\ref{eq:2.4}).
\end{tw}
The proof of Theorem \ref{tw:3.2} uses the same techniques and the same arguments as is used in the proof of existence theorem in three dimensional case, hence it will be omitted. The next Section will only very briefly present the main steps of the proof of Theorem \ref{tw:3.2}. The goal of this article is to prove a higher regularity of the displacement vector $u$ and the microrotation tensor $A$, which are the solutions of the system (\ref{eq:2.1}) in the sense of Definition \ref{de:3.1}. Let us denote by $C^{0,\alpha}([0,T];C^{0,\alpha}(\bar{\Omega};\R^3)$ the space of all H\"older continuous functions up to the boundary with exponent $\alpha>0$. The following theorem is the main result of this article. 
\begin{tw}
$\mathrm{(Main\; result)}$\\
\label{tw:3.3}
Let us  assume that for all $T>0$ the given data $f$, $g_D$, $A_D$ have the regularity
\begin{eqnarray*}
f\in W^{1,\infty}(0,T;L^2(\Omega;\R^3))\,,\quad g_D\in W^{1,\infty}(0,T;H^{\frac{1}{2}}(\partial\Omega;\R^3))\,,
\end{eqnarray*}
\begin{eqnarray*}
A_D\in H^{1}(0,T;H^{\frac{3}{2}}(\partial\Omega;\so(3)))
\end{eqnarray*}
and that there exists function $w\in W^{1,\infty}(0,T;H^1(\Omega;\R^3))$ such that $w_{t_{|_{\partial\Omega}}}=g_{D,t_{|_{\partial\Omega}}}$, satisfying 
$$\int\nolimits_{B(x_0,R)\cap\Omega}|\nabla\,w_t(x,t)|^2\,\di x\leq K\,R^{\gamma}\quad \mathrm{ for \;almost\; all\;} t\in (0,T)\,,$$
where $\gamma>0$ is any positive number and $B(x_0,R)\subset\R^2$ denotes the open ball with the center $x_0\in\R^2$ and the radius $R>0$ (the constant $K>0$ does not depend on the radius $R$). Additionally let us suppose that the initial data $(\ve^{p,0},b^0)\in L^2(\Omega;\SS)\times L^2(\Omega;\SS)$\\ satisfy (\ref{eq:3.3}) and (\ref{eq:3.4}). Then 
$$u\in C^{0,\alpha}([0,T];C^{0,\alpha}(\bar{\Omega};\R^3))\quad \mathrm{{\it and}}\quad A\in C^{0,\alpha}([0,T];C^{0,\alpha}(\bar{\Omega};\so(3)))$$ for $0<\alpha<1$, where the displacement vector $u$ and the microrotation tensor $A$ are the solutions (in the sense of Definition \ref{de:3.1}) of the system (\ref{eq:2.1})
with boundary conditions (\ref{eq:2.2}), (\ref{eq:2.3}) and initial condition (\ref{eq:2.4}).
\end{tw}
Theorem \ref{tw:3.3} presents the  first regularity result for non-monotone models from elasto-plasticity. The proof of Theorem \ref{tw:3.3} is based on the method of Morrey, which was presented in \cite{23} and \cite{24}. It is divided into three sections. First, we use the Yosida Approximation to the maximal monotone part of the inelastic constitutive equation and we very shortly show the main steps of the proof of Theorem \ref{tw:3.2}. Next, we prove the tube-filling condition (interior case). It will be the main part of the proof of Theorem \ref{tw:3.3}. Finally we show a Morrey's condition for the displacement vector up to the boundary $\partial\Omega$.
 
%%%%%%%%%%%%%%%%%%%%%%%%%%%%%%%%%%%%%%%%%%%%%%%%%%%%% Section 4 %%%%%%%%%%%%%%%%%%%%%%%%%%%%%%%%%%%%%%%%%%%%%%%%%%%%%%%%%%%%%%%%%%%%%%%%

\section{Existence for each Yosida approximation step}
We apply the Yosida approximation for the monotone part of the flow rule from (\ref{eq:2.1}) in order to get only a Lipschitz-nonlinearity in equation $(\ref{eq:2.1})_4$. We consider the following initial-boundary value problem
\renewcommand{\theequation}{\thesection.\arabic{equation}}
\setcounter{equation}{0}%
\begin{eqnarray}
\label{eq:4.1}
\mathrm{div}_x\, T^{\nu}&=&-f\,,\nn\\[1ex]
T^{\nu}&=&2\,\mu(\ve(u^{\nu})-\ve^{p,\nu})+2\,\mu_c(\mathrm{skew}\,(\nabla_x u^{\nu})-A^{\nu}) +\lambda\mathrm{tr}\,(\ve(u^{\nu})-\ve^{p,\nu})\,\id\,,\nn\\[1ex]
-l_c\,\Delta_x\,\mathrm{axl}\,(A^{\nu})&=&\mu_c\,\mathrm{axl}\,(\mathrm{skew}(\nabla_x u^{\nu})-A^{\nu})\,,\nn\\[1ex]
\ve^{p,\nu}_{t}&=&\frac{1}{\nu}\,\{|\dev\,(T_E^{\nu})-b^{\nu}|-\KK\}_{+}\frac{\dev\,(T_E^{\nu})-b^{\nu}}{|\dev\,(T_E^{\nu}) -b^{\nu}|}\,,\\[1ex]
T_E^{\nu}&=&2\,\mu(\ve(u^{\nu})-\ve^{p,\nu})+\lambda\mathrm{tr}\,(\ve(u^{\nu})-\ve^{p,\nu})\id,\nn\\[1ex]
b_t^{\nu}&=&c\,\ve^{p,\nu}_t-d\,|\ve_t^{p,\nu}|\,b^{\nu}.\nn
\end{eqnarray}
The above equations are studied for $x\in\Omega\subset \R^2$ and $t\in (0,T)$. $\nu>0$ is the Yosida approximation parameter and $\{\rho\}_{+}=\max\,\{0,\rho\}$, where $\rho$ is a 
scalar function.\\
The system (\ref{eq:4.1}) is considered with boundary conditions:
\begin{eqnarray}
\label{eq:4.2}
u^{\nu}(x,t)&=&g_D(x,t)\qquad \textrm{ for}\quad x\in\partial\Omega \quad\textrm{and}\quad t\geq 0\textrm{,}\nn\\[1ex]
A^{\nu}(x,t)&=&A_D(x,t)\qquad \textrm{ for}\quad x\in \partial \Omega \quad\textrm{and}\quad  t\geq 0
\end{eqnarray}
and initial conditions
\begin{eqnarray}
\label{eq:4.3}
\ve^{p,\nu}(x,0)=\ve^{p,0}(x),\qquad b^{\nu}(x,0)=b^0(x).
\end{eqnarray}
Denote by $\E^{\nu}(t)$ the total energy associated with the system (\ref{eq:4.1})
\begin{eqnarray}
\label{eq:4.4}
\E^{\nu}(u^{\nu}\ve^{\nu},\ve^{p,\nu},A^{\nu},b^{\nu})(t)&=&\int\nolimits_{\Omega}
\rho\,\psi^{\nu}\Big(u^{\nu}(x,t),\ve^{\nu}(x,t),\ve^{p,\nu}(x,t),A^{\nu}(x,t),b^{\nu}(x,t)\Big)\,\di x\nn\\[1ex]
&=&\int\nolimits_{\Omega}\Big(\mu\,\|\ve(u^{\nu})-\ve^{p,\nu}\|^2+
\mu_c\,\|\mathrm{skew}\,(\nabla_x u^{\nu})-A^{\nu}\|^2\\[1ex]
&+&\frac{\lambda}{2}\Big(\mathrm{tr}\,(\ve(u^{\nu})-\ve^{p,\nu})\Big)^2+2\,l_c\,\|\nabla_x\,\mathrm{axl}\,(A^{\nu})\|^2+\frac{1}{2c}\,\|b^{\nu}\|^2\Big)\,\di x\nn\,.
\end{eqnarray}
\begin{de}
\label{de:4.1}
Fix $T>0$. We say that a vector $(u,A,T,\ve^p,b)\in W^{1,\infty}(0,T;H^1(\Omega;\R^3)\times H^2(\Omega;\mathfrak{so}(3))\times L^2(\Omega;\R^9)\times L^2(\Omega;\SS)\times L^2(\Omega;\SS))$ is an $L^2$- strong solution of the system (\ref{eq:4.1}) if
\begin{enumerate}
\item $|\ve^p_t|\,b\in L^{\infty}(0,T;L^2(\Omega,\SS))$,
\item $|\dev\,\Big(2\mu(\ve(u(x,t))-\ve^p(x,t))\Big)-b(x,t)|\leq \KK$ for almost all $(x,t)\in \Omega\times (0,T)$,
\item the equations (\ref{eq:4.1}) are satisfied for almost all $(x,t)\in \Omega\times (0,T)$.
\end{enumerate}
\end{de}
The following lemma implies the $L^{\infty}$- boundedness of the backstress $b^{\nu}$.
\begin{lem}
\label{lem:4.2}
Fix $T>0$. Assume that $(u^{\nu},T^{\nu},A^{\nu},\ve^{p,\nu},b^{\nu})$ is an $L^2$- strong solution of the problem (\ref{eq:4.1}) and $|b^0(x)|\leq \frac{c}{d}$ for almost all $x\in\Omega$. Then for all $\nu>0$
$$|b^{\nu}(x,t)|\leq \frac{c}{d}\quad\mathrm{for\; a.\;e.\;}\quad (x,t)\in\Omega\times (0,T).$$ 
\end{lem}
For the proof of Lemma \ref{lem:4.2} we refer to \cite{7}. Now we propose the existence of solutions for each approximation step.
\begin{tw}
\label{tw:4.3}
Fix $T>0$. Suppose that all hypotheses of Theorem \ref{tw:3.2} are satisfied. Then for all $\nu>0$ there exists a unique $L^2$- strong solution (in the sense of Definition \ref{de:4.1})  
$$(u^{\nu},T^{\nu},A^{\nu},\ve^{p,\nu},b^{\nu})\in W^{1,\infty}\Big(0,T;H^1(\Omega;\R^3)\times L^2(\Omega;\R^9)\times H^2(\Omega;\mathfrak{so}(3))\times (L^2(\Omega;\SS))^2\Big)$$ satisfying the system (\ref{eq:4.1}) with boundary conditions (\ref{eq:4.2}) and initial condition (\ref{eq:4.3}).   
\end{tw}
The proof of Theorem \ref{tw:4.3} is the same as for tree-dimensional cases. It uses the same techniques as for the related Armstrong-Frederick model without Cosserat effects: see Section 4 of \cite{7}, therefore it will be omitted. For the complete proof of Theorem \ref{tw:4.3} we refer to \cite{12}.\\[1ex]
A fundamental tool in the proof of Theorem \ref{tw:3.2} is the following property of the energy function which results from our Cosserat modification:
\begin{tw}
\label{tw:4.4}
$\mathrm{(coerciveness\; of\; the\; energy)}$\\
(a) $\mathit{(the\; case\; with\; zero\; boundary\; data)}$\\[1ex]
For all $\nu>0$ the energy function (\ref{eq:4.4}) is elastically coercive with respect to $\nabla u$. 
This means that $\exists\; C_E>0$, $\forall\; u\in H^1_0(\Omega)$, $\forall\; A\in H^1_0(\Omega)$, $\forall\; \ve^p\in L^2(\Omega)$, $\forall\; b\in L^2(\Omega)$
$$\E^{\nu}(u,\ve,\ve^{p},A,b)\geq C_E\Big(\|u\|^2_{H^1(\Omega)}+\|A\|^2_{H^1(\Omega)}+\|b\|^2_{L^2(\Omega)}\Big).$$ 
(b) $\mathit{(the\; case\; with\; non-zero\; boundary\; data)}$\\[1ex]
Moreover, $\exists\; C_E>0$, $\forall\; g_D,\;A_D\in H^{\frac{1}{2}}(\partial\Omega)$, 
$\exists\; C_D>0$, $\forall\; \ve^p\in L^2(\Omega)$, $\forall\; b\in L^2(\Omega)$, $\forall\; u\in H^1(\Omega)$, $\forall\; A\in H^1(\Omega)$ 
with boundary conditions $u_{|_{\partial\Omega}}=g_D$ and $A_{|_{\partial\Omega}}=A_D$ it holds that
$$\E^{\nu}(u,\ve,\ve^{p},A,b)+C_D\geq C_E\Big(\|u\|^2_{H^1(\Omega)}+\|A\|^2_{H^1(\Omega)}+\|b\|^2_{L^2(\Omega)}\Big).$$
\end{tw}
For the proof of Theorem \ref{tw:4.4} we refer to the Theorem $3.2$ of the article \cite{8}.\\[1ex]
To pass to the limit in the system (\ref{eq:4.1}) and obtain the solution in the sense of Definition \ref{de:3.1} we need estimates for the time derivatives of the sequence $(u^{\nu},A^{\nu},T^{\nu},\ve^{p,\nu},b^{\nu})$. The article \cite{12} yields that the following energy estimate is sufficient to pass to the limit with the Yosida approximation.
\begin{tw}$\mathrm{(Energy\; estimate)}$\\
\label{tw:4.5}
Assume that the given data and initial data satisfies (\ref{eq:3.1}) - (\ref{eq:3.4}). Then for all $t\in(0,T)$ the following estimate
\begin{eqnarray*}
&&\int\nolimits_{\Omega}\frac{1}{2\nu}\{|\dev\,(T_E^{\nu})(t)-b^{\nu}(t)|-\KK\}_{+}^2\,\di x+\int\nolimits^t_0\int\nolimits_{\Omega}\D^{-1}\,T_{E,t}^{\nu}(\tau)\,T_{E,t}^{\nu}(\tau)\,\di x\,\di\tau\\[1ex]
&+& 2\,\mu_c\int\nolimits^t_0\int\nolimits_{\Omega}|\mathrm{skew}\,(\nabla_x u_t^{\nu}(\tau))-A_t^{\nu}(\tau)|^2\,\di x\,\di\tau
+4\,l_c\int\nolimits^t_0\int\nolimits_{\Omega}|\nabla\,\mathrm{axl}\,(A_t^{\nu}(\tau))|^2\,\di x\,\di\tau
\leq C(T)
\end{eqnarray*}
holds and $C(T)$ does not depend on $\nu>0$ (it depends only on the given data and the domain $\Omega$).
\end{tw}
For the proof of Theorem \ref{tw:4.5} we refer to the Section 5 of the paper \cite{12}. Theorem \ref{tw:4.5} and the elastic constitutive equations $(\ref{eq:4.1})_6$ imply that for a subsequence (again denoted by $\nu$) we have 
\begin{eqnarray}
\label{eq:4.5}
u^{\nu}\rightharpoonup u\quad &\mathrm{in}& \quad L^{\infty}(0,T;H^1(\Omega;\R^3))\,,\nn\\[1ex]
u^{\nu}_t\rightharpoonup u_t\quad &\mathrm{in}& \quad L^{2}(0,T;H^1(\Omega;\R^3))\,,\nn\\[1ex]
A^{\nu}\rightharpoonup A\quad &\mathrm{in}& \quad L^{\infty}(0,T;H^2(\Omega;\so(3)))\,,\nn\\[1ex]
A^{\nu}_t\rightharpoonup A_t\quad &\mathrm{in}& \quad L^{2}(0,T;H^2(\Omega;\so(3)))\,,\\[1ex]
T^{\nu}\rightharpoonup T\quad &\mathrm{in}& \quad L^{\infty}(0,T;L^2(\Omega;\R^9))\,,\nn\\[1ex]
T^{\nu}_t\rightharpoonup T_t\quad &\mathrm{in}& \quad L^{2}(0,T;L^2(\Omega;\R^9))\,,\nn\\[1ex]
b^{\nu}\rightharpoonup b\quad &\mathrm{in}& \quad L^{\infty}(0,T;L^2(\Omega;\SS))\,,\nn\\[1ex]
b^{\nu}_t\rightharpoonup b_t\quad &\mathrm{in}& \quad L^{2}(0,T;L^2(\Omega;\SS))\,.\nn
\end{eqnarray}
The informations contained in (\ref{eq:4.5}) are enough to pass to the limit in the Yosida approximation and get the solution in the sense of Definition \ref{de:3.1}. The details may be found in \cite{12}.

%%%%%%%%%%%%%%%%%%%%%%%%%%%%%%%%%%%%%%%%%%%%%%%%%%%%% Section 5 %%%%%%%%%%%%%%%%%%%%%%%%%%%%%%%%%%%%%%%%%%%%%%%%%%%%%%%%%%%%%%%%%%%%%%%%

\section{H\"older continuity for displacement and microrotation (interior case)}
This section is the main part of the proof of Theorem \ref{tw:3.3}.\\[1ex]
Let us denote by $B_R=B_R(x_0)=B(x_0,R)\subset \R^2$ the open ball with center $x_0\in\R^2$ and radius $R$. Moreover, let $B_{2R}\bl B_{R}=B(x_0,2R)\bl B(x_0,R)$ be the open annulus with center $x_0\in\R^2$. First we formulate two lemmas that will be useful later on.
\begin{lem}
\label{lem:5.1}
For each $1\leq p<\infty$ there exists a constant $C$, depending only on $n$ and $p$, such that
$$\int\nolimits_{B_R}|v(y)-v(z)|^p\,\di y\leq C\,R^{n+p-1}\int\nolimits_{B_R}|\nabla\, v(y)|^p\,|y-z|^{1-n}\,\di y$$
for all $B_R\subset \R^n$, $v\in C^1(B_R;\R^3)$ and $z\in B_R$.
\end{lem}
For the proof of Lemma \ref{lem:5.1} we refer to Section 4.5.2 of \cite{14}.
\begin{lem}
\label{lem:5.2}
Let $B_{2R}\subset \R^2$ be an open ball and $v\in H^1(B_{2R};\R^3)$. Then, there exists a constant $\tilde{C}$, not depending on $v$, such that
$$\Big(\int\nolimits_{B_{2R}}|v(y)-c_R|^2\,\di y\Big)^{\frac{1}{2}}\leq \tilde{C}\,R\,\Big(\int\nolimits_{B_{2R}}|\nabla v(y)|^2\,\di y\Big)^{\frac{1}{2}}\,,$$
where $$c_R=\frac{1}{|B_{2R}\bl B_R|}\,\int\nolimits_{B_{2R}\bl B_R}v(y)\,\di y\,.$$
\end{lem}
{\bf\em Proof:}\hspace{2ex} Without loss of generality we may assume that $v\in C^1(B_{2R};\R^3)$ and obtain
\begin{eqnarray*}
\int\nolimits_{B_{2R}}|v(y)-c_R|^2\,\di y&=& \int\nolimits_{B(x,2R)}\Big|\frac{1}{|B_{2R}\bl B_{R}|}\,\int\nolimits_{B_{2R}\bl B_{R}}(v(y)-v(z))\,\di z\Big|^2\,\di y\\[2ex]
&\leq&\int\nolimits_{B_{2R}}\frac{1}{|B_{2R}\bl B_{R}|}\,\int\nolimits_{B_{2R}\bl B_{R}}|v(y)-v(z)|^2\,\di z\,\di y\\[2ex]
&\leq&\int\nolimits_{B_{2R}}\frac{1}{|B_{2R}\bl B_{R}|}\,\int\nolimits_{B_{2R}}|v(y)-v(z)|^2\,\di z\,\di y\\[2ex]
&\leq&\qquad (\mathrm{Lemma}\; \ref{lem:5.1})\qquad\leq\\[2ex]
&\leq&\int\nolimits_{B_{2R}}C\,R\,\int\nolimits_{B_{2R}}|\nabla\, v(z)|^2|z-y|^{-1}\,\di z\,\di y\\[2ex]
&\leq&\hat{C}\,R^2\,\int\nolimits_{B_{2R}}|\nabla\, v(z)|^2\,\di z.\hspace{35ex}\Box
\end{eqnarray*}
Notice that Lemma \ref{lem:5.2} is the Poincar\'e inequality with a special constant $c_R$ (see Section 4.5.2 of \cite{14}, where the Poincar\'e inequality is proven with another constant). To prove H\"older continuity of the displacement vector first we prove the following tube-filling condition. The idea is taken from the articles \cite{15} and \cite{17}.
\begin{tw}
\label{tw:5.3}
Assume that the given data and initial data satisfy all hypotheses of \\Theorem \ref{tw:3.3}. Then there exists constants $C,\,K\,>0$ that do not depend on $\nu>0$ such that the tube-filling condition
\begin{eqnarray*}
\int\nolimits_0^t\int\nolimits_{B_R}|\nabla\, u^{\nu}_t(\tau)|^2\,\di x\,\di \tau\leq C\,\int\nolimits_0^t\int\nolimits_{B_{2R}\bl B_R}|\nabla\, u^{\nu}_t(\tau)|^2\,\di x\,\di \tau + K\,R^{\gamma}
\end{eqnarray*}
is satisfied for all balls $B_R=B(x_0,R)\subset B_{2R}=B(x_0,2R)\subset \Omega\subset \R^2$, $t\in(0,T)$ and $\gamma>0$. Moreover the constants $C,\,K>0$ do not depend on the radius $R.$
\end{tw}
{\bf\em Proof:}\hspace{2ex} Let us define the following function
\begin{displaymath}
\xi_0(s) = \left\{ \begin{array}{ll}
1 & s\in [-R,R],\\ 
\frac{2R-|s|}{R} & s\in [-2R,2R]\backslash [-R,R],\\
0 & s\notin [-2R,2R]
\end{array} \right.
\end{displaymath}
and let $\xi$ denotes the cutoff function defined by $\xi(x)=\xi_0(|x_0-x|)$. Compute the time derivative 
\renewcommand{\theequation}{\thesection.\arabic{equation}}
\setcounter{equation}{0}%
\begin{eqnarray} 
\label{eq:5.1} 
&&\frac{\di}{\di t}\,\Big (\int\nolimits_{B_{2R}}\xi^2\,\frac{1}{2\nu}\{|\dev\,(T_E^{\nu}(t))-b^{\nu}(t)|-\KK\}_{+}^2\,\di x\Big )\nn\\[1ex]
&=& \int\nolimits_{B_{2R}} \xi^2\,\ve^{p,\nu}_t(t)\,(\dev\,(T_{E,t}^{\nu}(t))-b^{\nu}_t(t))\,\di x
=(\textrm{by}\;\textrm{the}\;\textrm{elastic}\;\textrm{constitutive}\;\textrm{relation})\nn\\[1ex]
&=&\int\nolimits_{B_{2R}} \xi^2\,\ve^{\nu}_t(t)\,T_{E,t}^{\nu}(t)\,\di x-
\int\nolimits_{B_{2R}} \xi^2\,\D^{-1}\,T_{E,t}^{\nu}(t)\,T_{E,t}^{\nu}(t)\,\di x\nn\\[1ex]
&-&\int\nolimits_{B_{2R}} \xi^2\,\ve^{p,\nu}_t(t)\,b^{\nu}_t(t)\,\di x.
\end{eqnarray}
From Lemma \ref{lem:4.2} we conclude the following inequality
\begin{eqnarray}
\label{eq:5.2}  
\int\nolimits_{B_{2R}} \xi^2\,\ve^{p,\nu}_t(t)\,b^{\nu}_t(t)\,\di x&=& 
(\textrm{by}\;\textrm{the}\;\textrm{equation}\;\textrm{for}\;\textrm{the}\;\textrm{backstress})=\nn\\
&=&c\int\nolimits_{B_{2R}}\xi^2\, |\ve^{p,\nu}_t(t)|^2\,\di x
-d\int\nolimits_{B_{2R}}\xi^2\, |\ve^{p,\nu}_t(t)|\,\ve^{p,\nu}_t(t)\,b^{\nu}(t)\,\di x\nn\\[1ex]
&\geq& \int\nolimits_{B_{2R}}\xi^2\, |\ve^{p,\nu}_t(t)|^2\,(c-d\,|b^{\nu}(t)|)\,\di x\geq 0\,.
\end{eqnarray}
Notice that
\begin{eqnarray} 
\label{eq:5.3} 
\int\nolimits_{B_{2R}} \xi^2\,\ve^{\nu}_t(t)\,T_{E,t}^{\nu}(t)\,\di x&=&\int\nolimits_{B_{2R}} \xi^2\,\nabla\, u^{\nu}_t(t)\,T_{t}^{\nu}(t)\,\di x\nn\\[1ex]
&-&
2\,\mu_c\int\nolimits_{B_{2R}}\xi^2\,\Big(\mathrm{skew}\,(\nabla_x u_t^{\nu}(t))-A_t^{\nu}(t)\Big)\,\mathrm{skew}\,(\nabla_x u^{\nu}_t(t))\,\di x\nn\\[1ex]
&=&\int\nolimits_{B_{2R}} \xi^2\,\nabla\, u^{\nu}_t(t)\,T_{t}^{\nu}(t)\,\di x\nn\\[1ex]
&-&
2\,\mu_c\int\nolimits_{B_{2R}}\xi^2\,|\mathrm{skew}\,(\nabla_x u_t^{\nu}(t))-A_t^{\nu}(t)|^2\,\di x\nn\\[1ex]
&-&2\,\mu_c\int\nolimits_{B_{2R}}\xi^2\,\Big(\mathrm{skew}\,(\nabla_x u_t^{\nu}(t))-A_t^{\nu}(t)\Big)\,A_t^{\nu}(t)\,\di x.
\end{eqnarray} 
Integrating (\ref{eq:5.3}) with respect to time over $0$ to $t$ we have
\begin{eqnarray} 
\label{eq:5.4} 
\int\nolimits_0^t\int\nolimits_{B_{2R}}\xi^2\, \ve^{\nu}_t(\tau)\,T_{t}^{\nu}(\tau)\,\di x\,\di\tau&=& \int\nolimits_0^t\int\nolimits_{B_{2R}} \xi^2\,\nabla\, u^{\nu}_t(\tau)\,T_{t}^{\nu}(\tau)\,\di x\,\di\tau\nn\\[1ex]
&-&2\,\mu_c\int\nolimits_0^t\int\nolimits_{B_{2R}}\xi^2\,|\mathrm{skew}\,(\nabla_x u_t^{\nu}(\tau))-A_t^{\nu}(\tau)|^2\,\di x\,\di\tau\\[1ex]
&-&2\,\mu_c\int\nolimits_0^t\int\nolimits_{B_{2R}}\xi^2\,\Big(\mathrm{skew}\,(\nabla_x u_t^{\nu}(\tau))-A_t^{\nu}(\tau)\Big)\,A_t^{\nu}(\tau)\,\di x\,\di\tau\,.\nn
\end{eqnarray}
Let $\bar{u}_t^{\nu}$ be the average of $u_t^{\nu}$ on the set $B_{2R}\backslash B_R$. Integrating by parts in the first term on the right hand side of (\ref{eq:5.4}) we obtain
\begin{eqnarray} 
\label{eq:5.5} 
\int\nolimits_0^t\int\nolimits_{B_{2R}}\xi^2\, \nabla\, u^{\nu}_t(\tau)\,T_{t}^{\nu}(\tau)\,\di x\,\di\tau&=& \int\nolimits_0^t\int\nolimits_{B_{2R}} \xi^2\,(u^{\nu}_t(\tau)-\bar{u}_t^{\nu}(\tau))\,f_t(\tau)\,\di x\,\di\tau\\[1ex]
&-&2\int\nolimits_0^t\int\nolimits_{B_{2R}\bl B_R}\xi\, T_{t}^{\nu}(\tau)\cdot(u^{\nu}_t(\tau)-\bar{u}_t^{\nu}(\tau))\otimes\nabla\xi\,\di x\,\di\tau\nn\,.
\end{eqnarray}
The first term on the right hand side of (\ref{eq:5.5}) is estimated as follows
\begin{eqnarray} 
\label{eq:5.6} 
&&\int\nolimits_0^t\int\nolimits_{B_{2R}} \xi^2\,(u^{\nu}_t(\tau)-\bar{u}_t^{\nu}(\tau))\,f_t(\tau)\,\di x\,\di\tau \leq \int\nolimits_0^t\|u^{\nu}_t(\tau)-\bar{u}_t^{\nu}(\tau)\|_{L^2}\,\|\xi^2\, f_t(\tau)\|_{L^2}\,\di\tau\nn\\
&\leq&\;(\textrm{Lemma }\ref{lem:5.2})\;\leq \tilde{C}\,R\int\nolimits_0^t\|\nabla\, u^{\nu}_t(\tau)\|_{L^2}\,\di\tau
\leq\;(\textrm{Theorem }\ref{tw:4.5})\leq C(T)\,R
\end{eqnarray}
and the constants $\tilde{C}$, $C(T)>0$ do not depend on $\nu>0$. Using the Cauchy inequality with a small weight and applying Poincar$\mathrm{\acute{e}}$'s inequality to the second term on the right hand side of (\ref{eq:5.5}) we get 
\begin{eqnarray} 
\label{eq:5.7} 
&&\Bigg|\,2\int\nolimits_0^t\int\nolimits_{B_{2R}\backslash B_R}\xi\, T_{t}^{\nu}(\tau)\cdot(u^{\nu}_t(\tau)-\bar{u}_t^{\nu}(\tau))\otimes\nabla\,\xi \,\di x\,\di\tau\,\Bigg|\leq \alpha\int\limits_0^t\|\xi\, T_{t}^{\nu}(\tau)\|^2_{L^2}\,\di \tau\nn\\[1ex]
&+&\frac{\hat{C}(\alpha)}{R^2}\int\nolimits_0^t\int\nolimits_{B_{2R}\bl B_R}|u^{\nu}_t(\tau)-\bar{u}_t^{\nu}(\tau)|^2\,\di x\,\di \tau\nn\\[1ex]
&\leq&\alpha\int\nolimits_0^t\|\xi\, T_{t}^{\nu}(t)\|^2_{L^2}\,\di\tau+ C(\alpha)\int\nolimits_0^t\int\nolimits_{B_{2R}\bl B_R}|\nabla\,u^{\nu}_t(\tau)|^2\,\di x\,\di \tau\, ,
\end{eqnarray}
where $\alpha>0$ is any positive constant and the constants $\hat{C}(\alpha)$, $C(\alpha)>0$ do not depend on $\nu>0$. The last term on the right hand side of (\ref{eq:5.4}) is estimated as follows  
\begin{eqnarray} 
\label{eq:5.8} 
&&\Bigg|\,2\,\mu_c\int\nolimits_0^t\int\nolimits_{B_{2R}}\xi^2\,\Big(\mathrm{skew}\,(\nabla_x u_t^{\nu}(\tau))-A_t^{\nu}(\tau)\Big)\,A_t^{\nu}(t)\,\di x\,\di \tau\,\Bigg|\nn\\[1ex]
&\leq& \tilde{\alpha}\int\nolimits_0^t\int\nolimits_{B_{2R}}\xi^2\,|\mathrm{skew}\,(\nabla_x u_t^{\nu}(\tau))-A_t^{\nu}(\tau)|^2\,\di x\,\di t\nn\\[1ex]
&+& C(\tilde{\alpha})\int\nolimits_0^t\int\nolimits_{B_{2R}}\xi^2\,|A_t^{\nu}(\tau)|^2\,\di x\,\di \tau\, ,
\end{eqnarray}
where $\tilde{\alpha}>0$ is any positive constant and the constant $C(\tilde{\alpha})>0$ does not depend on $\nu>0$. Notice that the sequence $\{A^\nu_t\}_{\nu>0}$ is bounded in $L^{2}(0,T;H^2(\Omega;\so(3)))$, then from the Rellich - Kondrachov Theorem (cf. \cite{1}) with $n=2$ we obtain that the sequence $\{A^\nu_t\}_{\nu>0}$ is bounded in $L^{2}(0,T;C^{0,\beta}(\Omega;\so(3)))$, where $0<\beta<1$. Hence 
\begin{eqnarray} 
\label{eq:5.9}
\int\nolimits_0^t\int\nolimits_{B_{2R}}\xi^2\,|A_t^{\nu}(\tau)|^2\,\di x\,\di\tau\leq C\,R^2\int\nolimits_0^T\max_{x\in B_{2R}}|A_t^{\nu}(t)|^2\,\di t\leq \tilde{C}\,R^2\,.
\end{eqnarray}
Integrating (\ref{eq:5.1}) with respect to time and using (\ref{eq:5.2}) - (\ref{eq:5.9}) we obtain the following inequality
\begin{eqnarray}
\label{eq:5.10}
&&\int\nolimits_{B_{2R}}\xi^2\,\frac{1}{2\nu}\{|\dev\,(T_E^{\nu}(t))-b^{\nu}(t)|-\KK\}_{+}^2\,\di x+\int\nolimits^t_0\int\nolimits_{B_{2R}}\xi^2\,\D^{-1}\,T_{E,t}^{\nu}(\tau)\,T_{E,t}^{\nu}(\tau)\,\di x\,\di\tau\nn\\[1ex]
&+& 2\,\mu_c\int\nolimits^t_0\int\nolimits_{B_{2R}}\xi^2\,|\mathrm{skew}\,(\nabla_x u_t^{\nu}(\tau))-A_t^{\nu}(\tau)|^2\,\di x\,\di\tau
\,\,\leq\,\, \alpha\int\nolimits_0^t\int\nolimits_{B_{2R}}\xi^2\,|T_{E,t}^{\nu}(\tau)|^2\,\di x\,\di\tau\nn\\[1ex]
&+& \alpha\int\nolimits_0^t\int\nolimits_{B_{2R}}\xi^2\,|\mathrm{skew}\,(\nabla_x u_t^{\nu}(\tau))-A_t^{\nu}(\tau)|^2\,\di x\,\di\tau + C(\alpha)\int\nolimits_0^t\int\nolimits_{B_{2R}\bl B_R}|\nabla\, u^{\nu}_t(\tau)|^2\,\di x\,\di\tau\nn\\[1ex]
&+&C\,R^{\gamma}
+\int\nolimits_{B_{2R}}\xi^2\,\frac{1}{2\nu}\{|\dev\,(T_E^{\nu}(0))-b^{\nu}(0)|-\KK\}_{+}^2\,\di x\, ,
\end{eqnarray}
where the constants $C(\alpha),\,C>0$ do not depend on $\nu>0$ and the radius $R$ ($\gamma$ is any positive constant). Observe that $T^{\nu}(0)\in L^2(\Omega;\S)$ is the unique solution of the problem
\begin{eqnarray}
\label{eq:5.35}
\mathrm{div}_x\, T^{\nu}(x,0)&=&-f(x,0)\,,\nn\\[1ex]
-l_c\,\Delta_x\,\mathrm{axl}\,(A^{\nu}(x,0))&=&\mu_c\,\mathrm{axl}\,(\mathrm{skew}\,(\nabla_x u^{\nu}(x,0))-A^{\nu}(x,0))\,,\\[1ex]
u^{\nu}(x,0)_{|_{\partial\Omega}}=g_D(x,0)\,, && A^{x,0}(x,0)_{|_{\partial\Omega}}=A_D(x,0)\,,\nn
\end{eqnarray}
where 
\begin{eqnarray*}
T^{\nu}(x,0)&=&2\mu\,(\ve(u^{\nu}(x,0))-\ve^{p,0}(x))+2\mu_c(\mathrm{skew}\,(\nabla_x u^{\nu}(x,0))-A^{\nu}(x,0))\\[1ex]
&+&\lambda\,\mathrm{tr}(\ve(u^{\nu}(x,0))-\ve^{p,0}(x))\,\id\,,
\end{eqnarray*}
which implies that $\dev\,(T^{\nu}_E(0))=\dev\,(T^0_E)$ and $b^{\nu}(0)=b^0$. From the assumption (\ref{eq:3.3}) we conclude that the last term on the right hand side of (\ref{eq:5.10}) is equal to zero.  Choosing in (\ref{eq:5.10}) $\alpha>0$ sufficiently small we arrive that 
\begin{eqnarray}
\label{eq:5.11}
&&\int\nolimits_{B_{2R}}\xi^2\,\frac{1}{2\nu}\{|\dev\,(T_E^{\nu}(t))-b^{\nu}(t)|-\KK\}_{+}^2\,\di x+\int\nolimits^t_0\int\nolimits_{B_{2R}}\xi^2\,\D^{-1}\,T_{E,t}^{\nu}(\tau)\,T_{E,t}^{\nu}(\tau)\,\di x\,\di\tau\nn\\[1ex]
&+& 2\,\mu_c\int\nolimits^t_0\int\nolimits_{B_{2R}}\xi^2\,|\mathrm{skew}\,(\nabla_x u_t^{\nu}(\tau))-A_t^{\nu}(\tau)|^2\,\di x\,\di \tau\nn\\[1ex]
&\leq& C\int\nolimits^t_0\int\nolimits_{B_{2R}\bl B_R}|\nabla\, u^{\nu}_t(\tau)|^2\,\di x\,\di\tau + C\,R^{\gamma}\,.
\end{eqnarray}
The inequality (\ref{eq:5.11}) implies the following inequality
\begin{eqnarray}
\label{eq:5.12}
&&\int\nolimits^t_0\int\nolimits_{B_{2R}}\xi^2\,\D^{-1}\,T_{E,t}^{\nu}(\tau)\,T_{E,t}^{\nu}(\tau)\,\di x\,\di\tau
+ 2\,\mu_c\int\nolimits^t_0\int\nolimits_{B_{2R}}\xi^2\,|\mathrm{skew}\,(\nabla_x u_t^{\nu}(\tau))-A_t^{\nu}(\tau)|^2\,\di x\,\di\tau\nn\\[1ex]
&\leq&  C\int\nolimits^t_0\int\nolimits_{B_{2R}\bl B_R}|\nabla\, u^{\nu}_t(\tau)|^2\,\di x\,\di\tau+C\,R^{\gamma}.
\end{eqnarray}
Notice also that
\begin{eqnarray}
\label{eq:5.13}
&&2\,\mu_c\int\nolimits^t_0\int\nolimits_{B_{2R}}\xi^2\,|\mathrm{skew}\,(\nabla_x u_t^{\nu}(\tau))-A_t^{\nu}(\tau)|^2\,\di x\,\di\tau\nn\\[1ex]
&=&2\,\mu_c\int\nolimits^t_0\int\nolimits_{B_{2R}}\xi^2\,|\mathrm{skew}\,(\nabla_x u_t^{\nu}(\tau))|^2\,\di x\,\di\tau
+2\,\mu_c\int\nolimits^t_0\int\nolimits_{B_{2R}}\xi^2\,|A_t^{\nu}(\tau)|^2\,\di x\,\di\tau\nn\\[1ex]
&-&4\,\mu_c\int\nolimits^t_0\int\nolimits_{B_{2R}}\xi\,\mathrm{skew}\,(\nabla_x u_t^{\nu}(\tau))\cdot\xi\, A_t^{\nu}(\tau)\,\di x\,\di\tau\nn\\[1ex]
&\geq & \mu_c\int\nolimits^t_0\int\nolimits_{B_{2R}}\xi^2\,|\mathrm{skew}\,(\nabla_x u_t^{\nu}(\tau))|^2\,\di x\,\di\tau - 
2\,\mu_c\int\nolimits^t_0\int\nolimits_{B_{2R}}\xi^2\,|A_t^{\nu}(\tau)|^2\,\di x\,\di\tau.\qquad\qquad
\end{eqnarray}
From the observation $\mathrm{div}\,u^{\nu}_t=\mathrm{tr}\,\ve(u^{\nu}_t)-\mathrm{tr}\,\ve^{p,\nu}_t$ we have
\begin{eqnarray}
\label{eq:5.14}
\int\nolimits^t_0\int\nolimits_{B_{2R}}\xi^2\,|\mathrm{div}\, u_t^{\nu}(\tau)|^2\,\di x\,\di\tau\leq \int\nolimits^t_0\int\nolimits_{B_{2R}}\xi^2\,|T_{E,t}^{\nu}(\tau)|^2\,\di x\,\di t\, ,
\end{eqnarray}
From the assumptions on the elasticity tensor $\D$ we know that
\begin{eqnarray}
\label{eq:5.15}
\int\nolimits^t_0\int\nolimits_{B_{2R}}\xi^2\,\D^{-1}\,T_{E,t}^{\nu}(\tau)\,T_{E,t}^{\nu}(\tau)\,\di x\,\di\tau\geq D\int\nolimits^t_0\int\nolimits_{B_{2R}}\xi^2\,|T_{E,t}^{\nu}(\tau)|^2\,\di x\,\di\tau\,.
\end{eqnarray}
and the constant $D>0$ does not depend on $\nu>0$. Using the expressions (\ref{eq:5.13}) - (\ref{eq:5.15}) in (\ref{eq:5.12}) we get 
\begin{eqnarray}
\label{eq:5.16}
&&D\int\nolimits^t_0\int\nolimits_{B_{2R}}\xi^2\,|\mathrm{div}\, u_t^{\nu}(\tau)|^2\,\di x\,\di\tau+\mu_c\int\nolimits^t_0\int\nolimits_{B_{2R}}\xi^2\,|\mathrm{skew}\,(\nabla_x u_t^{\nu}(\tau))|^2\,\di x\,\di\tau\nn\\[1ex]
&\leq& C\int\nolimits^t_0\int\nolimits_{B_{2R}\bl B_R}|\nabla\, u^{\nu}_t(\tau)|^2\,\di x\,\di\tau + C\,R^{\gamma}+
2\,\mu_c\int\nolimits^t_0\int\nolimits_{B_{2R}}\xi^2\,|A_t^{\nu}(\tau)|^2\,\di x\,\di\tau\,.\quad
\end{eqnarray}
The last term on the right hand side of (\ref{eq:5.16}) is estimated in the same way as in (\ref{eq:5.9}). To complete the proof we need to estimate the expression
$$\int\nolimits^t_0\int\nolimits_{B_{2R}}\xi^2\,|\mathrm{div}\, u_t^{\nu}(\tau)|^2\,\di x\,\di \tau+\int\nolimits^t_0\int\nolimits_{B_{2R}}\xi^2\,|\mathrm{curl}\, u_t^{\nu}(\tau))|^2\,\di x\,\di\tau\,.$$
Let us denote by 
$$c_R=\frac{1}{|B_{2R}\bl B_{R}|}\,\int\limits_{B_{2R}\bl B_{R}}u^{\nu}_t(x,t)\,\di x\,,$$
then
\begin{eqnarray}
\label{eq:5.17}
\mathrm{curl}\,(\xi\,(u_t^{\nu}-c_R))&=&\nabla\,\xi\times(u_t^{\nu}-c_R)+\xi\,\mathrm{curl}\,u^{\nu}_t\,,\nn\\[1ex]
\mathrm{div}\,(\xi\,(u-c_R))&=&\nabla\,\xi\cdot(u_t^{\nu}-c_R)+\xi\,\mathrm{div}\,u_t^{\nu}
\end{eqnarray}
and
\begin{eqnarray}
\label{eq:5.18}
|\xi\,\mathrm{curl}\,u_t^{\nu}|^2&=&|\mathrm{curl}\,(\xi\,(u_t^{\nu}-c_R))|^2+|\nabla\,\xi\times(u_t^{\nu}-c_R)|^2\nn\\[1ex] &-&2\,\mathrm{curl}\,(\xi\,(u_t^{\nu}-c_R))\cdot\nabla\,\xi\times(u_t^{\nu}-c_R)\nn\\[1ex]
&\geq&|\mathrm{curl}\,(\xi\,(u_t^{\nu}-c_R))|^2+|\nabla\,\xi\times(u_t^{\nu}-c_R)|^2\nn\\[1ex]
&-&2\,|\mathrm{curl}\,(\xi\,(u_t^{\nu}-c_R))|\,|\nabla\,\xi\times(u_t^{\nu}-c_R)|\nn\\[1ex]
&\geq&|\mathrm{curl}\,(\xi\,(u_t^{\nu}-c_R))|^2+|\nabla\,\xi\times(u_t^{\nu}-c_R)|^2- \epsilon\,|\mathrm{curl}\,(\xi\,(u_t^{\nu}-c_R))|^2\nn\\[1ex]
&-&\frac{C}{\epsilon}\,|\nabla\,\xi\times(u_t^{\nu}-c_R)|^2\,\geq\,\mathrm{(for\, sufficiently\, small\, epsilon)\,}\,\geq\nn\\[1ex]
&\geq&C\,|\mathrm{curl}\,(\xi\,(u_t^{\nu}-c_R))|^2- \tilde{C}\,|\nabla\,\xi\times(u_t^{\nu}-c_R)|^2\,.
\end{eqnarray}
In the same manner as in (\ref{eq:5.18}) we arrive at the following inequality\\
\begin{eqnarray}
\label{eq:5.19}
|\xi\,\mathrm{div}\,u_t^{\nu}|^2\geq C\,|\mathrm{div}\,(\xi\,(u-c_R))|^2-
\tilde{C}\,|\nabla\,\xi\cdot(u_t^{\nu}-c_R)|^2\, .\\[1ex]\nn
\end{eqnarray}
Using (\ref{eq:5.17}) - (\ref{eq:5.19}) in (\ref{eq:5.16}) we have 
\begin{eqnarray*}
&&C_1\int\nolimits^t_0\int\nolimits_{B_{2R}}\Big(|\mathrm{div}\,(\xi\,(u_t^{\nu}(\tau)-c_R))|^2+ |\mathrm{curl}\,(\xi\,(u_t^{\nu}(\tau)-c_R))|^2\Big)\,\di x\,\di\tau\\[1ex]
&\leq&C\int\nolimits^t_0\int\nolimits_{B_{2R}\bl B_R}|\nabla\, u^{\nu}_t(\tau)|^2\,\di x\,\di\tau  + C\,R^{\gamma}+
C_2\,\int\nolimits^t_0\int\nolimits_{B_{2R}\bl B_{R}}|\nabla\,\xi|^2\,|u_t^{\nu}(\tau)-c_R|^2\,\di x\,\di\tau\,,
\end{eqnarray*}
where the constants $C$, $\tilde{C}$, $C_1$ and $C_2$ do not depend on $\nu>0$. Notice that the function $\xi\,(u_t^{\nu}-c_R)\in H^1_0(B_{2R};\R^3)$ for almost all $t>0$, hence the well-known estimate [\cite{18}, p.36]\\
$$\|\nabla u\|^2_{L^2}\leq C^{\mathrm{curl}}_{\mathrm{div}}\Big(\|\mathrm{div}\; u\|^2_{L^2}+
\|\mathrm{curl}\; u\|^2_{L^2}\Big)\quad \mathrm{for\, all}\quad u\in H^1_0\\[1ex]$$
(the constant $C^{\mathrm{curl}}_{\mathrm{div}}$ does not depend on $u$) implies the following inequality
\begin{eqnarray}
\label{eq:5.20}
&&C_1\int\nolimits^t_0\int\nolimits_{B_{2R}}|\nabla\,(\xi\,(u_t^{\nu}(\tau)-c_R))|^2\,\di x\,\di\tau
\leq C\int\nolimits^t_0\int\nolimits_{B_{2R}\bl B_R}|\nabla\, u^{\nu}_t(\tau)|^2\,\di x\,\di\tau\nn\\[1ex]
&+&C_2\,\int\nolimits^t_0\int\nolimits_{B_{2R}\bl B_{R}}|\nabla\,\xi|^2\,|u_t^{\nu}(\tau)-c_R|^2\,\di 
x\,\di+C\,R^{\gamma}.
\end{eqnarray}
The expression
\begin{eqnarray*}
|\nabla\,(\xi\,(u_t^{\nu}-c_R))|^2=|\nabla\,\xi\otimes(u_t^{\nu}-c_R)|^2+|\xi\,\nabla\,u_t^{\nu}|^2
+2\,\nabla\,\xi\otimes(u_t^{\nu}-c_R)\cdot\xi\,\nabla\,u_t^{\nu}
\end{eqnarray*}
yields the inequality
\begin{eqnarray}
\label{eq:5.21}
&&C_1\int\nolimits^t_0\int\nolimits_{B_{2R}}\xi^2\,|\nabla\,u_t^{\nu}(\tau)|^2\,\di x\,\di\tau
\leq C\int\nolimits^t_0\int\nolimits_{B_{2R}\bl B_R}|\nabla u^{\nu}_t(\tau)|^2\,\di x\,\di\tau + C\,R^{\gamma}\nn\\[1ex]
&+& C_2\,\int\nolimits^t_0\int\nolimits_{B_{2R}\bl B_{R}}|\nabla\,\xi|^2\,|u_t^{\nu}(\tau)-c_R|^2\,\di x\,\di\tau\nn\\[1ex]
&-&2\,\int\nolimits^t_0\int\nolimits_{B_{2R}}\nabla\,\xi\otimes(u_t^{\nu}(\tau)-c_R)\cdot\xi\,\nabla\,u_t^{\nu}(\tau)\,\di x\,\di\tau\nn\\[1ex]
&\leq& C\int\nolimits^t_0\int\nolimits_{B_{2R}\bl B_R}|\nabla\, u^{\nu}_t(\tau)|^2\,\di x\,\di\tau+C\,R^{\gamma}\nn\\[1ex]
&+&C_2\,\int\nolimits^t_0\int\nolimits_{B_{2R}\bl B_{R}}|\nabla\,\xi|^2\,|u_t^{\nu}(\tau)-c_R|^2\,\di x\,\di\tau\nn\\
&+& C(a)\int\nolimits^t_0\int\nolimits_{B_{2R}\bl B_R}|\nabla\,\xi|^2\,|u_t^{\nu}(\tau)-c_R|^2\,\di x\,\di\tau + a\,\int\nolimits^t_0\int\nolimits_{B_{2R}}\xi^2\,|\nabla\,u_t^{\nu}(\tau)|^2\,\di x\,\di\tau\qquad\;
\end{eqnarray}
for all $a>0$. Choosing in (\ref{eq:5.21}) $a>0$ sufficiently small we obtain
\begin{eqnarray}
\label{eq:5.22}
C_1\int\nolimits^t_0\int\nolimits_{B_{2R}}\xi^2\,|\nabla\,u_t^{\nu}(\tau)|^2\,\di x\,\di\tau
&\leq& C\int\nolimits^t_0\int\nolimits_{B_{2R}\bl B_{R}}|\nabla\,u_t^{\nu}(\tau)|^2\,\di x\,\di\tau + C\,R^{\gamma}\nn\\[1ex]
&+&C_3\,\int\nolimits^t_0\int\nolimits_{B_{2R}\bl B_{R}}|\nabla\,\xi|^2\,|u_t^{\nu}(\tau)-c_R|^2\,\di x\,\di\tau\,.\qquad
\end{eqnarray}
The last term on the right hand side of (\ref{eq:5.22}) is estimated using Poincar\'e's inequality
\begin{eqnarray}
\label{eq:5.23}
\int\nolimits^t_0\int\nolimits_{B_{2R}\bl B_{R}}|\nabla\,\xi|^2\,|u_t^{\nu}(\tau)-c_R|^2\,\di x\,\di\tau &\leq& \frac{C}{R^2}\,\int\nolimits^t_0\int\nolimits_{B_{2R}\bl B_{R}}|u_t^{\nu}(\tau)-c_R|^2\,\di x\,\di\tau\nn\\[1ex]
&\leq&C\,\int\nolimits^t_0\int\nolimits_{B_{2R}\bl B_{R}}|\nabla\,u_t^{\nu}(\tau)|^2\,\di x\,\di\tau\, ,
\end{eqnarray}
where the constant $C>0$ does not depend on $\nu>0$ and the radius $R>0$. Applying (\ref{eq:5.23}) in (\ref{eq:5.22}) and the fact that $\xi\equiv 1$ on $B_R$ we complete the proof.$\mbox{}$ \hfill $\Box$\\[2ex]
To prove local H\"older continuity for the displacement vector $u$ we will show the following Morrey's condition (see for instance \cite{6}): in the case of two dimensions it is in the form  
$$\int\nolimits_0^t\int\nolimits_{B_{R}}|\nabla\,u_t^{\nu}(\tau)|\,\di x\,\di\tau\leq K\,R^{1+\gamma}\qquad \forall\quad B_R\subset\Omega\subset\R^2\, ,$$ 
where $0<\gamma<1$, $t\in(0,T)$ and the constant $K>0$ does not depend on the radius $R>0$. To obtain the above inequality we use the Widman's hole filling trick from the articles \cite{27} and \cite{15}.
\begin{tw}
\label{tw:5.4}
Suppose that all hypotheses of Theorem \ref{tw:3.3} holds. Then there exists $\alpha\in (0,1)$ such that
\begin{eqnarray*}
\int\nolimits^t_0\int\nolimits_{B_R}|\nabla\, u^{\nu}_t(\tau)|^2\,\di x\,\di\tau\leq 2^{2\alpha}\,\frac{R^{2\alpha}}{R_0^{2\alpha}}\,\Big(\int\nolimits^t_0\int\nolimits_{B_{R_0}}|\nabla\, u^{\nu}_t(\tau)|^2\,\di x\,\di\tau + \tilde{K}\,R_0^{\gamma}\Big)
\end{eqnarray*}
for all balls $B_R=B_R(x_0,R)\subset\Omega\subset \R^2$ and $t\in(0,T)$, where $x_0\in\Omega'\Subset\Omega$ and $2\,R\leq R_0=\frac{1}{2}\,\mathrm{dist}\,(\Omega',\partial\Omega)$. The constants $\tilde{K}$ (independent on the radius $R>0$) and $R_0$ do not depend on $\nu>0$ and $\gamma>0$ is any positive constant. 
\end{tw}
{\bf\em Proof:}\hspace{2ex} From Theorem \ref{tw:5.3} we conclude the following tube filling condition
\begin{eqnarray}
\label{eq:5.24}
\int\nolimits^t_0\int\nolimits_{B_R}|\nabla\, u^{\nu}_t(\tau)|^2\,\di x\,\di\tau\leq C\int\nolimits^t_0\int\nolimits_{B_{2R}\bl B_R}|\nabla\, u^{\nu}_t(\tau)|^2\,\di x\,\di\tau + K\,R^{\gamma}\,,
\end{eqnarray}
where $\gamma>0$ and the constant $K$ does not depend on $\nu$. We add to both sides of (\ref{eq:5.24}) the expression (filling the hole) $C\int\nolimits^t_0\int\nolimits_{B_R}|\nabla\, u^{\nu}_t(\tau)|^2\,\di x\,\di\tau$ and we obtain
\begin{eqnarray}
\label{eq:5.25}
\int\nolimits^t_0\int\nolimits_{B_R}|\nabla\, u^{\nu}_t(\tau)|^2\,\di x\,\di\tau\leq \frac{C}{1+C}\,\int\nolimits^t_0\int\nolimits_{B_{2R}}|\nabla\, u^{\nu}_t(\tau)|^2\,\di x\,\di\tau + K\,R^{\gamma}.
\end{eqnarray}
For $j\geq 1$ we set $R_j=R_0\,2^{-j}$ and by iteration we deduce that 
\begin{eqnarray}
\label{eq:5.26}
&&\int\nolimits^t_0\int\nolimits_{B_{R_N}}|\nabla\, u^{\nu}_t(\tau)|^2\,\di x\,\di\tau\leq \frac{C}{1+C}\,\int\nolimits^t_0\int\nolimits_{B_{R_{2N}}}|\nabla\, u^{\nu}_t(\tau)|^2\,\di x\,\di\tau + K\,R^{\gamma}_N\nn\\[1ex]
&\leq&\frac{C}{1+C}\,\Big(\frac{C}{1+C}\,\int\nolimits^t_0\int\nolimits_{B_{R_{4N}}}|\nabla\, u^{\nu}_t(\tau)|^2\,\di x\,\di\tau + K\,R^{\gamma}_{2N}\Big)+ K\,R^{\gamma}_N\leq\;\;...\;\;\leq\nn\\[1ex]
&\leq&\Big(\frac{C}{1+C}\Big)^N\,\int\nolimits^t_0\int\nolimits_{B_{R_{0}}}|\nabla\, u^{\nu}_t(\tau)|^2\,\di x\,\di\tau+K\,\sum\limits^N_{k=1}R_k^{\gamma}\,\Big(\frac{C}{1+C}\Big)^{N-k}.
\end{eqnarray}
Let us choose $\alpha$ such that $\frac{C}{1+C}=2^{-2\alpha}$, then we obtain that 
\begin{eqnarray}
\label{eq:5.27}
K\,\sum\limits^N_{k=1}R_0^{\gamma}\,2^{-\gamma\, k}\,2^{-2\,\alpha(N-k)}=K\,R_0^{\gamma}\,2^{-2\,\alpha\, N}\,\sum\limits^N_{k=1}2^{-(\gamma - 2\alpha)k}\leq \tilde{K}\,R_0^{\gamma}\,(2^{-N})^{2\alpha}
\end{eqnarray}
for $\alpha< \frac{\gamma}{2}$ and $\tilde{K}<\infty$. Using (\ref{eq:5.27}) in (\ref{eq:5.26}) we obtain
\begin{eqnarray}
\label{eq:5.28}
\int\nolimits^t_0\int\nolimits_{B_{R_N}}|\nabla\, u^{\nu}_t(\tau)|^2\,\di x\,\di\tau\leq (2^{-N})^{2\alpha}\,\Big(\int\nolimits^t_0\int\nolimits_{B_{R_0}}|\nabla\, u^{\nu}_t(\tau)|^2\,\di x\,\di\tau + \tilde{K}\,R_0^{\gamma}\Big)\,.
\end{eqnarray}
From the assumption $R\leq R_0/2$ it is possible to find $j\geq 1$ such that $R_{j+1}\leq R\leq R_{j}$ and $R_j\leq 2R$. Hence $(2R)^{-1}\leq (\frac{R_0}{2^j})^{-1}$ and
\begin{eqnarray*}
R^{-2\alpha}=(2R)^{-2\alpha}\,2^{2\alpha}\leq 2^{2\alpha}\,\Big(\frac{R_0}{2^j}\Big)^{-2\alpha}\quad&/&\cdot\, (2^j)^{2\alpha}\\[2ex]
(2^j)^{2\alpha}\,R^{-2\alpha}\leq 2^{2\alpha}\,R_0^{-2\alpha}\quad &/&\cdot\, R^{2\alpha}\\[2ex]
(2^j)^{2\alpha}\leq 2^{2\alpha}\,\frac{R^{2\alpha}}{R_0^{2\alpha}}\,.\quad &&
\end{eqnarray*}
The above inequality finishes the proof.$\mbox{}$ \hfill $\Box$\\[2ex]
Theorem \ref{tw:5.4} yields that the velocity of the displacement vector satisfies Morrey's condition for $n=2$ but we are unable to prove H\"older continuity of the velocity. We show it for the displacement - the idea was taken from the article \cite{17}. From the Theorem \ref{tw:4.5} and coerciveness of the total energy we know that $\nabla\, u_t^{\nu}\in H^1(0,T;L^2(\Omega;\R^9))$, hence  
$$\nabla\, u^{\nu}(x,t)=\nabla\, u(x,0)+\int\nolimits_0^t \nabla\,u_t^{\nu}(x,\tau)\,\di\tau\quad \textrm{for all}\quad 0\leq t\leq T$$
and
$$\int\nolimits_{B_R}|\nabla\, u^{\nu}(x,t)|\,\di x\leq\int\nolimits_{B_R}|\nabla\, u(x,0)|\,\di x+\int\nolimits_0^t \int\nolimits_{B_R}|\nabla\,u_t^{\nu}(x,\tau)|\,\di x\,\di\tau\,,$$
where the ball $B_R\subset\Omega$ is the same as in the Theorem \ref{tw:5.4}. Observe that the function $u(x,0)$ is the unique solution of the elliptic system (\ref{eq:5.35}) and it does not depend on $\nu>0$. The general regularity theory for elliptic systems (developed by K. O Widman in \cite{27} and C. B Morrey in \cite{24}) implies that $\int\nolimits_{B_R}|\nabla\, u(x,0)|\,\di x\leq \tilde{K}\,R^{1+\tilde{\alpha}}$ for some $0<\tilde{\alpha}<1$ and the constant $\tilde{K}>0$ does not depend on the radius $R$. Using Theorem \ref{tw:5.4} we infer that the function $u^{\nu}$ satisfies Morrey's condition i.e. 
$$\int\nolimits_{B_R}|\nabla\, u^{\nu}(x,t)|\,\di x\leq K\,R^{1+\alpha},$$
where the constants $K,\alpha$ and the ball $B_R$ come from Theorem \ref{tw:5.4}. Hence\\ $u^{\nu}\in L^{\infty}(0,T;C^{0,\alpha}(\Omega',\R^3))$. The weak convergence of the sequence $\{\nabla u^{\nu}\}_{\nu>0}$ in $L^{\infty}(L^2)$ implies that the function $u$ satisfies Morrey's condition and $u\in L^{\infty}(0,T;C^{0,\alpha}(\Omega',\R^3))$.  To prove H\"older continuity with respect to time we estimate the following difference (the idea is taken again from \cite{17})
\begin{eqnarray*}
|u(x_0,t_1)-u(x_0,t_2)|&\leq& \Big|u(x_0,t_1)-\fint\nolimits_{B_R}u(x,t_1)\,\di x\Big|+ \Big|u(x_0,t_2)-\fint\nolimits_{B_R}u(x,t_2)\,\di x\Big|\nn\\
&+& \Big|\fint\nolimits_{B_R}u(x,t_1)\,\di x-\fint\nolimits_{B_R}u(x,t_2)\,\di x\Big|\, .
\end{eqnarray*}
The H\"older continuity with respect to space of the function $u$ implies that 
\begin{eqnarray*}
 \Big|\,u(x_0,t_i)-\fint\nolimits_{B_R}u(x,t_i)\,\di x\,\Big|\leq K\,R^{\alpha}\quad (i=1,2)
\end{eqnarray*}
for almost all $t>0$. Let us choose $R>0$ such that $|t_2-t_1|^{\frac{1}{2}}=R^{\alpha+1}$, then we obtain
\begin{eqnarray*}
|u(x_0,t_1)-u(x_0,t_2)|&\leq& 2\,K\,R^{\alpha}+\frac{C}{R^2}\,\Big|\int\nolimits_{t_1}^{t_2}\int\nolimits_{B_R}u_t(x,t)\,\di x\,\di t\Big|\nn\\[1ex]
&\leq& 2\,K\,R^{\alpha}+ C\,R^{-1}\,\int\nolimits_{t_1}^{t_2}\Big(\int\nolimits_{B_R}|u_t(x,t)|^2\,\di x\Big)^{\frac{1}{2}}\,\di t\nn\\[1ex]
&\leq& 2\,K\,R^{\alpha}+ C\,R^{-1}\,|t_2-t_1|^{\frac{1}{2}}\,\Big(\int\nolimits_{t_1}^{t_2}\int\nolimits_{B_R}|u_t(x,t)|^2\,\di x\,\di t\Big)^{\frac{1}{2}}\,.\qquad\;
\end{eqnarray*}
We know that $u_t\in L^2(0,T;H^1(\Omega;\R^3))$, hence the choice of the radius $R$ yields that
\begin{eqnarray*}
|u(x_0,t_1)-u(x_0,t_2)|\leq K\,|t_2-t_1|^{\frac{\alpha}{2(\alpha+1)}}\,.
\end{eqnarray*}
In the same manner we can prove the H\"older continuity with respect to time of the function $A\in L^{\infty}(0,T;H^2(\Omega;\so(3)))$ because from the Sobolev imbedding Theorem (cf. \cite{1}) we have that $A\in L^{\infty}(0,T;C^{0,\alpha}(\bar{\Omega};\so(3)))$ and $A_t\in L^{2}(0,T;H^2(\Omega;\so(3)))$.

%%%%%%%%%%%%%%%%%%%%%%%%%%%%%%%%%%%%%%%%%%%%%%%%%%%%% Section 6 %%%%%%%%%%%%%%%%%%%%%%%%%%%%%%%%%%%%%%%%%%%%%%%%%%%%%%%%%%%%%%%%%%%%%%%%

\section{H\"older continuity of the displacement up to the boundary}
To prove H\"older continuity up to the boundary we need some assumptions about the boundary $\partial\Omega$. At the beginning of the Section 2 we established that $\partial\Omega$ is Lipschitz but we can also assume that $\partial\Omega$ satisfies the so called " Winer type condition ", cf. \cite{15}.\\[3ex]
{\large\bf\em Proof of Theorem \ref{tw:3.3}:}\hspace{2ex} The proof is divided into two steps:\\[1ex]
{\bf Step 1:}\hspace{2ex} Assume that $B_{2R}=B(x_0,2R)=B_{2R}(x_0)\subset\Omega\subset\R^2$. If we consider the same cutoff function as in the proof of Theorem \ref{tw:5.3}, then the values of u on the boundary are irrelevant and the proof is the same as in the last Section.\\[1ex]
{\bf Step 2:}\hspace{2ex} Assume that $B_{2R}\cap\partial\Omega\neq\emptyset$. Let us take a point $x'_0\in\,\partial\Omega$ such that $x'_0\in\partial\Omega\cap B_{2R}(x_0)$, where $B_R(x_0)\subset B_{4R}(x'_0)$. Let us consider a standard Lipschitz continuous cutoff function $\tau$ such that 
\renewcommand{\theequation}{\thesection.\arabic{equation}}
\setcounter{equation}{0}%
\begin{eqnarray*}
\tau&=&1\;\;\mathrm{on}\;\; B_{4R}(x_0')\,,\nn\\[1ex]
\tau&=&0\;\;\mathrm{on}\;\; \R^2\,\bl\, B_{8R}(x_0')\,,\\[1ex]
|\nabla\tau|&\leq& C\, R^{-1}\;\;\mathrm{on}\;\; B_{8R}(x_0')\,\bl\,B_{4R}(x_0') \,.\nn
\end{eqnarray*}
Regularity with respect to time of the data $g_{D}$ implies that there exists function $w\in W^{1,\infty}(0,T;H^1(\Omega;\R^3))$ such that $w_{t_{|_{\partial\Omega}}}=g_{D,t}$. Differentiating with respect to time the equations $(\ref{eq:4.1})_1$ and $(\ref{eq:4.1})_3$, next multiplying equation $(\ref{eq:4.1})_1$ by $\tau^2\,(u^{\nu}_t-w_t)$, equation $(\ref{eq:4.1})_3$ by $\tau^2\,\mathrm{axl}\,A^{\nu}_t$ and integrating those two equations with respect to space we obtain
\begin{eqnarray}
\label{eq:6.1}
\int\nolimits_{B_{8R}(x_0')\cap\Omega}\tau^2\,\mathrm{div}_x\, T^{\nu}_t\,(u^{\nu}_t-w_t)\,\di x&=& -\int\nolimits_{B_{8R}(x_0')\cap\Omega}\tau^2\,f_t\,(u^{\nu}_t-w_t)\,\di x\,,\nn\\[1ex]
-l_c\,\int\nolimits_{B_{8R}(x_0')\cap\Omega}\tau^2\,\Delta_x\,\mathrm{axl}\,(A_t^{\nu})\,\mathrm{axl}\,(A_t^{\nu})\,\di x&=&\\[1ex]
\mu_c\,\int\nolimits_{B_{8R}(x_0')\cap\Omega}&\tau^2&\mathrm{axl}\,(\mathrm{skew}(\nabla_x u^{\nu})-A^{\nu})\,\mathrm{axl}\,(A_t^{\nu})\,\di x\, .\nn
\end{eqnarray}
Integrating by parts in the first equation of (\ref{eq:6.1}) we have 
\begin{eqnarray}
\label{eq:6.2}
\int\nolimits_{B_{8R}(x_0')\cap\Omega}\tau^2\, T^{\nu}_t\,(\ve(u^{\nu}_t)-\ve(w_t))\,\di x&=& -\int\nolimits_{B_{8R}(x_0')\cap\Omega}\tau^2\,f_t\,(u^{\nu}_t-w_t)\,\di x\nn\\[1ex]
&-&2\int\nolimits_{T_R}\tau\, T_{t}^{\nu}\cdot (u^{\nu}_t-w_t)\otimes\nabla\tau\,\di x\, ,\nn\\[1ex]
-l_c\,\int\nolimits_{B_{8R}(x_0')\cap\Omega}\tau^2\,\Delta_x\,\mathrm{axl}\,(A_t^{\nu})\,\mathrm{axl}\,(A_t^{\nu})\,\di x&=&\\[1ex]
\mu_c\,\int\nolimits_{B_{8R}(x_0')\cap\Omega}&\tau^2&\mathrm{axl}\,(\mathrm{skew}(\nabla_x u^{\nu})-A^{\nu})\,\mathrm{axl}\,(A_t^{\nu})\,\di x\, ,\nn
\end{eqnarray}
where $T_R:=B_{8R}(x_0')\cap\Omega\,\bl\,B_{4R}(x_0')$. Adding those two equations we get
\begin{eqnarray}
\label{eq:6.3}
&&2\,\mu \int\nolimits_{B_{8R}(x_0')\cap\Omega}\tau^2\, |\ve(u^{\nu}_t)-\ve^{p,\nu}_t|^2\,\di x +\lambda\int\nolimits_{B_{8R}(x_0')\cap\Omega}\tau^2\, \Big(\mathrm{tr}\,(\ve(u^{\nu}_t)-\ve^{p,\nu}_t)\Big)^2\,\di x\nn\\[1ex]
&+&\int\nolimits_{B_{8R}(x_0')\cap\Omega}\tau^2\ve^{p,\nu}_t\,T^{\nu}_{E,t}\,\di x + 2\mu_c\int\nolimits_{B_{8R}(x_0')\cap\Omega}\tau^2\,|\mathrm{skew}(\nabla_x u_t^{\nu}(t))-A_t^{\nu}(t)|^2\,\di x\nn\\[1ex]
&=& -\int\nolimits_{B_{8R}(x_0')\cap\Omega}\tau^2\,f_t\,(u^{\nu}_t-w_t)\,\di x
-2\int\nolimits_{T_R}\tau\, T_{t}^{\nu}\cdot (u^{\nu}_t-w_t)\otimes\nabla\tau\,\di x\, ,\\[1ex]
&&+\int\nolimits_{B_{8R}(x_0')\cap\Omega}\tau^2\,T^{\nu}_t\, \ve(w_t)\,\di x-l_c\,\int\nolimits_{B_{8R}(x_0')\cap\Omega}\tau^2\,\Delta_x\,\mathrm{axl}\,(A_t^{\nu})\,\mathrm{axl}\,(A_t^{\nu})\,\di x\, .\nn
\end{eqnarray}
Integrating (\ref{eq:6.3}) with respect to time we obtain
\begin{eqnarray}
\label{eq:6.4}
&&\int\nolimits_{B_{8R}(x_0')\cap\Omega}\tau^2\,\frac{1}{2\nu}\{|\dev\,(T_E^{\nu}(t))-b^{\nu}(t)|-\KK\}_{+}^2\,\di x\nn\\[1ex]
&+&2\mu \int\nolimits_0^t\int\nolimits_{B_{8R}(x_0')\cap\Omega}\tau^2\, |\ve(u^{\nu}_t)-\ve^{p,\nu}_t|^2\,\di x\,\di\tau + \lambda\int\nolimits_0^t\int\nolimits_{B_{8R}(x_0')\cap\Omega}\tau^2\, \Big(\mathrm{tr}\,(\ve(u^{\nu}_t)-\ve^{p,\nu}_t)\Big)^2\,\di x\,\di\tau\nn\\[1ex]
&+&2\mu_c\int\nolimits_0^t\int\nolimits_{B_{8R}(x_0')\cap\Omega}\tau^2\,|\mathrm{skew}\,(\nabla_x u_t^{\nu})-A_t^{\nu}|^2\,\di x\,\di\tau\nn\\[1ex]
&=& -\int\nolimits_0^t\int\nolimits_{B_{8R}(x_0')\cap\Omega}\tau^2\,f_t\,(u^{\nu}_t-w_t)\,\di x\,\di\tau
-2\int\nolimits_0^t\int\nolimits_{T_R}\tau\, T_{t}^{\nu}\cdot (u^{\nu}_t-w_t)\otimes\nabla\tau\,\di x\,\di\tau \\
&+&\int\nolimits_0^t\int\nolimits_{B_{8R}(x_0')\cap\Omega}\tau^2\,T^{\nu}_t\, \ve(w_t)\,\di x\,\di\tau -l_c\,\int\nolimits_0^t\int\nolimits_{B_{8R}(x_0')\cap\Omega}\tau^2\,\Delta_x\,\mathrm{axl}\,(A_t^{\nu})\,\mathrm{axl}\,(A_t^{\nu})\,\di x\,\di\tau\nn\\[1ex]
&+&\int\nolimits_{B_{8R}(x_0')\cap\Omega}\tau^2\,\frac{1}{2\nu}\{|\dev\,(T_E^{\nu}(0))-b^{\nu}(0)|-\KK\}_{+}^2\,\di x\, .\nn
\end{eqnarray}
The same argument as in the proof of Theorem \ref{tw:5.3} yields that the last term on the right hand side of (\ref{eq:6.4}) equals zero. The first term on the right hand side of (\ref{eq:6.4}) is estimated as follows
\begin{eqnarray}
\label{eq:6.5} 
\int\nolimits_0^t\int\nolimits_{B_{8R}(x_0')\cap\Omega} \tau^2\,(u^{\nu}_t-w_t)\,f_t\,\di x\,\di\tau &\leq& \int\nolimits_0^t\|u^{\nu}_t-w_t\|_{L^2}\|\tau^2\,f_t(t)\|_{L^2}\,\di\tau\nn\\
&\leq& \tilde{C}R\,\int\nolimits_0^t\|\nabla\, u^{\nu}_t-\nabla\, w_t\|_{L^2}\,\di\tau
\leq C(T)\,R\, .\qquad
\end{eqnarray}
Here we apply the standard Poincar$\mathrm{\acute{e}}$ inequality, because $u_t^{\nu}-w_t=0$ on the set\\ $B_{8R}(x_0')\cap\partial\Omega$ which has positive measure ($\partial\Omega$ is Lipschitz hence it satisfies the sphere condition - see for example Section 1.1.3 of \cite{6}). Using Cauchy's inequality with a small weight and applying Poincar$\mathrm{\acute{e}}$'s inequality to the second term on the right hand side of (\ref{eq:6.4}) we conclude 
\begin{eqnarray}
\label{eq:6.6} 
2\,\int\nolimits_0^t\int\nolimits_{T_R}\tau\, T_{t}^{\nu}\cdot (u^{\nu}_t-w_t)\otimes\nabla\xi\, \di x\,\di\tau&\leq& a\,\int\nolimits_0^t\|\tau\, T_{t}^{\nu}\|^2_{L^2}\,\di\tau
+\frac{\hat{C}(a)}{R^2}\,\int\nolimits_0^t\int\nolimits_{T_R}|u^{\nu}_t-w_t|^2\,\di x\,\di t\nn\\
&\leq& a\int\nolimits_0^t\|\tau\, T_{t}^{\nu}\|^2_{L^2}\,\di\tau\nn\\
&+& C(a)\int\nolimits_0^t\int\nolimits_{T_R}|\nabla u^{\nu}_t-\nabla w_t|^2\di x\,\di\tau\, ,
\end{eqnarray}
where $a>0$ is any positive constant. One before last integral on the right hand side of (\ref{eq:6.4}) is estimated in the same way as in the proof of Theorem \ref{tw:5.3} because $A^{\nu}$ is the $L^2$- strong solution of (\ref{eq:4.1}). The estimates (\ref{eq:6.5}), (\ref{eq:6.6}) and the Cauchy inequality with a small weigh used in the third term of (\ref{eq:6.4}) implies the following inequality 
\begin{eqnarray}
\label{eq:6.7} 
&&\int\nolimits_{B_{8R}(x_0')\cap\Omega}\tau^2\,\frac{1}{2\nu}\,\{|\dev\,(T_E^{\nu}(t))-b^{\nu}(t)|-\KK\}_{+}^2\,\di x\nn\\[1ex]
&+&2\mu \int\nolimits_0^t\int\nolimits_{B_{8R}(x_0')\cap\Omega}\tau^2\, |\ve(u^{\nu}_t)-\ve^{p,\nu}_t|^2\,\di x\,\di\tau +\lambda\int\nolimits_0^t\int\nolimits_{B_{8R}(x_0')\cap\Omega}\tau^2\, \Big(\mathrm{tr}\,(\ve(u^{\nu}_t)-\ve^{p,\nu}_t)\Big)^2\,\di x\,\di\tau\nn\\[1ex]
&+&2\mu_c\int\nolimits_0^t\int\nolimits_{B_{8R}(x_0')\cap\Omega}\tau^2\,|\mathrm{skew}\,(\nabla_x u_t^{\nu})-A_t^{\nu}|^2\,\di x\,\di\tau\nn\\[1ex]
&\leq& C\,\int\nolimits_0^t\int\nolimits_{T_R}|\nabla\, u^{\nu}_t(t)|^2\,\di x\,\di\tau + a\,\int\nolimits_0^t\|\tau\, T_{t}^{\nu}\|^2_{L^2}\,\di\tau\nn\\[1ex]
&+&\tilde{C}\,\int\nolimits_0^t\int\nolimits_{B_{8R}(x_0')\cap\Omega}|\nabla\, w_t|^2\di x\,\di\tau\ + K\,R^{\gamma}
\end{eqnarray}
where $a,\,\gamma>0$ are arbitrary positive constants and the constants $C,\,\tilde{C}>0$ do not depend on $\nu>0$. Choosing in (\ref{eq:6.7}) $a>0$ suitably small, using the additional assumption on the function $w_t$ and the calculations from the proof of Theorem \ref{tw:5.3} we arrive at the following tube filling condition
\begin{eqnarray}
\label{eq:6.8} 
\int\nolimits_0^t\int\nolimits_{B_{4R}(x_0')}|\nabla\, u^{\nu}_t|^2\,\di x\,\di\tau\leq C\int\nolimits_0^t\int\nolimits_{T_R}|\nabla\, u^{\nu}_t|^2\,\di x\,\di\tau + K\,R^{\gamma},
\end{eqnarray}
where $\gamma>0$ is some positive constant. Let us note that 
\begin{eqnarray}
\label{eq:6.9} 
\int\nolimits_0^t\int\nolimits_{B_{R}(x_0)\cap\Omega}|\nabla\, u^{\nu}_t|^2\,\di x\,\di\tau\leq\int\nolimits_0^t\int\nolimits_{B_{4R}(x_0')\cap\Omega}|\nabla\, u^{\nu}_t|^2\,\di x\,\di\tau
\end{eqnarray}
and
\begin{eqnarray}
\label{eq:6.10} 
\int\nolimits_0^t\int\nolimits_{T_R}|\nabla\, u^{\nu}_t|^2\,\di x\,\di\tau\leq \int\nolimits_0^t\int\nolimits_{B_{16R}(x_0)\cap\Omega\bl B_R(x_0)}|\nabla\, u^{\nu}_t|^2\,\di x\,\di\tau\,.
\end{eqnarray}
The inequalities (\ref{eq:6.9}) and (\ref{eq:6.10}) yield that 
\begin{eqnarray*}
\int\nolimits_0^t\int\nolimits_{B_{R}(x_0)\cap\Omega}|\nabla\, u^{\nu}_t|^2\di x\,\di\tau\leq\int\nolimits_0^t\int\nolimits _{B_{16R}(x_0)\cap\Omega\bl B_R(x_0)}|\nabla\, u^{\nu}_t|^2\di x\,\di\tau+K\,R^{\gamma}\, ,
\end{eqnarray*}
therefore we can apply the tube-filling trick with $R_j=R_0\,16^{-j}$ (in the same manner as in the proof of Theorem \ref{tw:5.4}) and get finally
\begin{eqnarray}
\int\nolimits_{B_{R}(x_0)\cap\Omega}|\nabla\, u^{\nu}|^2\di x\leq K\,R^{2\alpha}
\end{eqnarray}
for all balls $B(x_0,R)\subset \R^2$, where $K<\infty$ and $0<\alpha<1$. The "Morrey's Dirichlet growth theorem" (see for example \cite{19}) implies that $u\in C^{0,\alpha}([0,T];C^{0,\alpha}(\bar{\Omega};\R^3))$.$\mbox{}$ \hfill $\Box$\\[5ex]
{\bf Acknowledgments: }This work has been supported by the European Union in
the framework of the European Social Fund through the Warsaw University of Technology Development Programme.
\footnotesize{
}
\end{document}